\documentclass{article}
\usepackage{amsmath}
\usepackage{amsfonts}
\usepackage[top=1.5in,bottom=1.5in,left=1.5in,right=1.5in]{geometry}

\numberwithin{equation}{section}
\newtheorem{thm}{Theorem}[section]
\newtheorem{cor}[thm]{Corollary}
\newtheorem{prop}[thm]{Proposition}

\newtheorem{lem}[thm]{Lemma}

\newcommand{\be}{\begin{equation}}
\newcommand{\ee}{\end{equation}}
\newcommand{\ben}{\begin{enumerate}}
\newcommand{\een}{\end{enumerate}}

\newcommand{\xiaosihao}{\fontsize{10pt}{\baselineskip}\selectfont}

\newcommand{\qed}{\hspace*{\fill}Q.E.D.}  
\title{\Large  Finsler metrics of weakly isotropic flag curvature}
\author{Benling Li
\\
\xiaosihao\it
Department of Mathematics
, Ningbo University \\
\xiaosihao\it Ningbo, Zhejiang Province 315211, P.R. China\\
\xiaosihao\it libenling@nbu.edu.cn
}

\date{ Final version Oct. 13, 2017, to appear in CAG (Accepted)}

\begin{document}
\maketitle

\begin{abstract}
Finsler metrics of scalar flag curvature play an important role to show the complexity and richness of Finsler geometry.
In this paper, on an $n$-dimensional manifold $M$ we study the Finsler metric $F=F(x,y)$ of scalar flag curvature ${\bf K} = {\bf K}(x,y)$ and discover some equations ${\bf K}$ should be satisfied. As an application,  we mainly study the metric $F$ of weakly isotropic flag curvature. We prove that in this case, $F$ must be a Randers metric when $dim(M) \geq 3$ or be of constant flag curvature. Without the restriction on the dimension, the same result is obtained for projectively flat Finsler metrics of  weakly isotropic flag curvature.
\end{abstract}

\vskip 5mm
\textbf{MR Subject Classification (2000)}  53B40, 53C60.

\section{Introduction}

  The flag curvature of Finsler metrics is a natural analogue of sectional curvature of Riemannian metrics. Let
$F =F (x,y)$ be a Finsler metric on an $n$-dimensional manifold $M$. The flag curvature ${\bf K} ={\bf K}(\Pi, y)$ of $F$ is a function
of "flag" $\Pi \subset T_{x}M$ and "flag pole" $y \in T_{x}M$  at $x$ with $y \in \Pi$.  A Finsler metric is of {\it scalar flag curvature}
${\bf K}= {\bf K}(x,y)$ if the flag curvature is independent of the "flags" $\Pi \subset T_{x}M$ for each "flag pole" $y$ at $x$ with $y\in \Pi$.
As a special class in Finsler geometry, with the quadratic restriction, if a Riemannian metric is of scalar flag curvature, then the flag curvature must be independent of the "flag pole", i.e. ${\bf K} = {\bf K}(x)$. By Schur's Lemma, in dimension $n\geq 3$, ${\bf K} = constant$. The classification of  Riemannian metrics of constant flag curvature (sectional curvature) is well-known. However, the classification of  Finsler metrics of scalar flag curvature (even constant flag curvature) is far from being understood. Then a natural inverse problem arises: how to determine the Finsler metric $F = F(x,y)$ when its flag curvature ${\bf K} = {\bf K}(x,y)$ is given?

  In this paper, we first establish two equations of the scalar flag curvature ${\bf K}$ should be satisfied in Lemma \ref{lemK1} and \ref{lemK2}. From these two equations, one can see that there are some restrictions on ${\bf K}$. As an application,
  we study the Finsler metrics of {\it weakly isotropic flag curvature}
  \be \label{Kiso}
  {\bf K} = \frac{3\theta}{F} + \sigma, \ee
 where $\theta = \theta_i(x) y^i$ is a $1$-form and $\sigma=\sigma(x)$ is a scalar function.  It is easy to see that when $\theta =0$, then ${\bf K} = \sigma(x)$ and hence ${\bf K} = constant$ in dimension $n\geq3$ by the Schur's Lemma. The classification of Finsler metrics of constant flag curvature is still unknown though many metrics of constant flag curvature were found \cite{JDG}\cite{Bryant}\cite{Bry2}\cite{Li2}\cite{LiShen}\cite{MSY}\cite{ShenProj0}. When $\theta \neq 0$, it seems more general than the previous case. That's why the name "weakly isotropic" flag curvature is given.
 The flag curvature in (\ref{Kiso}) was first considered in \cite{CMS} when the authors studied Finsler metrics with isotropic S-curvature. They proved that if a Finsler metric with isotropic S-curvature and of scalar flag curvature ${\bf K}={\bf K}(x,y)$, then ${\bf K}$ must be in the form (\ref{Kiso}). In the past years, many Finsler metrics of weakly isotropic flag curvature are discovered. However, all of them are Randers metrics when the 1-form $\theta \neq 0$.
 In \cite{ChSh}, X. Cheng and Z. Shen classified Randers metrics  of weakly isotropic flag curvature when $dim(M) \geq3$. These Randers metrics can be expressed by
 \be \label{Randers} F = \frac{\sqrt{ ( 1 - \| W \|^2_h) h^2 + (W_i y^i)^2 }}{ 1 - \| W \|^2_h }  - \frac{W_i y^i}{1 - \| W \|^2_h}, \ee
 where
 \[ h = \frac{\sqrt{ |y|^2 + \mu (|x|^2|y|^2 - \langle x,y \rangle^2) }}{1+\mu |x|^2}, \]
\[ W = -2\Big[ (\delta \sqrt{1+\mu |x|^2} + \langle a,x\rangle)x - \frac{a|x|^2}{\sqrt{1+\mu |x|^2} +1}  \Big] + x Q + b + \mu\langle b,x\rangle x, \]
$\delta$, $\mu$ are constant, $Q = (q_j^{\ i})$ is an anti-symmetric matrix and $a$, $b \in R^n$ are constant vectors. In this case, the flag curvature
is given by
\[ {\bf K} = 3\frac{c_{x^m}y^m}{F} + \sigma,  \]
where
\[ c = \frac{\delta + \langle a,x\rangle}{\sqrt{1+\mu |x|^2}},  \ \ \ \sigma = \mu -c^2 - 2 c_{x^m}W^m. \]

 Then it is natural to ask that if a Finsler metric of weakly isotropic flag curvature ($\theta \neq 0$) must be a Randers metric or not?
In the following theorem, we give the positive answer to this question when the dimension of the manifold is not less than 3.
\begin{thm}\label{thm1}
Let $F$ be a Finsler metric on an $n$-dimensional manifold $M$ ($n\geq 3$) of weakly isotropic flag curvature. Then $F$ must be  a Randers metric in (\ref{Randers}) or be of constant flag curvature.
\end{thm}

When $\theta =0$, then ${\bf K} = \sigma(x)$ and by Schur's Lemma ${\bf K} = constant.$ When $\theta \neq 0$, the proof is based on Lemma \ref{lemK1} and Lemma \ref{lemK2} obtained in Section \ref{sectionLemma}.
We first prove the 1-form $\theta$ is closed and then find a quadratic equation $F$ should be satisfied.
By the above theorem and Lemma \ref{lemShen}, we immediately get the following corollary.
\begin{cor}
Let $F$ be a Finsler metric on an $n$-dimensional manifold $M$ ($n\geq 3$) of scalar flag curvature and having
 almost  isotropic S-curvature, i.e., \[ S = (n+1)(c F + \eta), \]
where $c = c(x)$ is a scalar function and $\eta = \eta_i(x)y^i$ is a closed 1-form. If $c(x)$ is not a constant,  then $F$ is a Randers metric in (\ref{Randers}).
\end{cor}

  If there is no restriction on flag curvature, some results of Randers metrics having almost S-curvature can be found in \cite{ChSh}\cite{Deng}\cite{Xing}.

  In Theorem \ref{thm1}, we need the condition that $dim(M)\geq 3 $. Then it is natural to ask how about the case
  when $dim(M) =2$? Till now, we still do not know how to discuss this case in general. However, as we known that
  all projectively flat Finsler metrics are of scalar flag curvature. For this special class we obtain the following result
  without the restriction on the dimension.
\begin{thm}\label{thm2}
Let $F$ be a  projectively flat Finsler metric of weakly isotropic flag curvature
\[ {\bf K} = \frac{3\theta}{F} + \sigma,\]
where $\theta = \theta_i(x) y^i \neq 0$ is a $1$-form and $\sigma=\sigma(x)$ is a scalar function. Then  $F$ is a Randers metric.
\end{thm}

It is known that
Randers metrics of weakly isotropic flag curvature must be with isotropic S-curvature.
 In \cite{CMS}, projectively flat Randers metric with isotropic S-curvature is totally determined.
When $\theta =0$, then ${\bf K} = \sigma(x)$. It is known that projectively flat Finsler metrics with  ${\bf K} = \sigma(x)$ must be of constant flag curvature.
In dimension greater than two, it follows from
the Schur's Lemma. In dimension two, it is proved by L. Berwald. In \cite{Li2}, we have given the classification of projectively flat Finsler metrics with constant flag curvature.

The following corollary is obvious by Theorem \ref{thm2} and Lemma \ref{lemShen}.
\begin{cor}
Let $F$ be a locally projectively flat Finsler metric with almost isotropic S-curvature, i.e., \[ S = (n+1)(c F + \eta), \]
where $c = c(x)$ is a scalar function and $\eta = \eta_i(x)y^i$ is a closed 1-form. If $c(x)$ is not a constant,  then $F$ is a Randers metric.
\end{cor}

Then we ask the following \\
{\bf Open Problem}: On a 2-dimensional manifold, is there any Finsler metric which is not a Randers metric  of weakly isotropic flag curvature in (\ref{Kiso}) ($\theta \neq 0$)?

\section{Preliminaries}
Let $M$ be a $C^{\infty}$  $n$-dimensional  manifold. A {\it Finlser metric} $F=F(x,y)$ on $M$ is a $C^{\infty}$ function
on $TM \setminus \{0\}$ with the following properties: (i) $F\geq0$ and $F(x,y)=0$ if and only if
 $y=0$; (ii) $F$ is a positively homogeneous function of degree one, i.e., $F (x, \lambda y) = \lambda F(x,y)$, $\lambda\geq0$;
 (iii) $F$ is strongly convex, i.e., for any $y \neq 0$, the matrix $g_{ij}:= \frac{1}{2} [F^2]_{y^iy^j}$
is positive definite. Specially, if $g_{ij} = g_{ij}(x)$, then $F$ is called a {\it Riemannian metric}, i.e. $F = \sqrt{g_{ij}(x)y^i y^j}$.
{\it Randers metrics} compose another group of simplest metrics in Finsler geometry, which arise from many
areas in mathematics, physics and biology \cite{Antonelli}. They can be expressed in the
 form $F =\alpha +\beta$, where $\alpha = \sqrt{a_{ij}(x)y^i y^j}$ is a Riemannian metric and $\beta = b_i(x)y^i$ is a 1-form with
$\| \beta \|_{\alpha} <1$ for any point.

   Consider a Finsler metric $F=F(x,y)$ on an open domain $\mathcal{U} \subset R^n$. The geodesics of $F$ are characterized by the following ODEs:
\[ \frac{d^2 x^i}{dt^2} + 2 G^i (x, \frac{dx}{dt}) = 0,\]
where $ G^i = G^i(x,y)$ are called {\it geodesic coefficients }
given by
\[ G^i = \frac{1}{4} g^{il} \Big\{  \frac{\partial^2 [F^2]}{\partial x^m \partial y^l} y^m -\frac{\partial [F^2]}{\partial x^l} \Big\}. \]

   The Riemann curvature ${\bf R}_y = R^i_{\ k} \frac{\partial}{\partial x^i} \otimes dx^k$ is defined by $G^i$ as the following
    \[  R^i_{\  k}:= 2\frac{\partial G^i}{\partial x^k} - \frac{\partial^2 G^i}{\partial x^l \partial y^k}y^l + 2 G^l \frac{\partial^2 G^i}{\partial y^l \partial y^k} - \frac{\partial G^i}{\partial y^l} \frac{\partial G^l}{\partial y^k}.  \]

    As an extension of sectional curvature in Riemann geometry, for each tangent plane $\Pi\subset T_xM$ and $y\in \Pi$, the
\emph{flag curvature} of $(\Pi,y)$ is defined by
\begin{eqnarray}
\nonumber
{\bf K}(\Pi,y)=\frac{g_{im}R^i_{\ k} u^ku^m}{F^2g_{ij}u^iu^j-[g_{ij}y^iu^j]^2},
\end{eqnarray}
where $\Pi=span\{y,u\}$.
 A Finsler metric $F$ is of {\it scalar flag
curvature} if its flag curvature ${\bf K}(\Pi,y)={\bf K}(x,y)$ is
independent of the tangent plane $\Pi$. In this case
\[ R^i_{\ k}={\bf K} \{F^2\delta^i_k-F \frac{\partial F}{\partial y^k} y^i\}.\]
If $F$ is a  Riemannian metric, the flag curvature ${\bf K}(\Pi,y)={\bf K}(\Pi)$ is
independent of $y$. Finsler metric $F$ is said to be of {\it
isotropic flag curvature} if ${\bf K}=K(x)$. If ${\bf K}$ is a constant, then $F$ is said to be of
{\it constant flag curvature}.

  In Finsler geometry, there are some non-Riemannian quantities. There quantities always vanish for Riemannian
  metrics. Such as Cartan torsion, mean Cartan torsion, Berwald curvature, Landsberg curvature,   S-curvature and
  etc.

The {\it Cartan torsion} ${\bf C}_y = C_{ijk} dx^i \otimes dx^j \otimes dx^k$ is defined by
\[ C_{ijk} = \frac{1}{2}\frac{\partial g_{ij}}{\partial y^k} = \frac{1}{4} \frac{\partial^3 F^2}{\partial y^i \partial y^j \partial y^k}. \]
And the {\it mean Cartan torsion} ${\bf I}_y = I_k dx^k$ is defined by
\[ I_k = g^{ij} C_{ijk}. \]
 It is obvious that $G^i = \frac{1}{2}\Gamma^i_{jk}(x)y^j y^k$ are
quadratic in $y$ for a Riemannian metric, where $\Gamma^i_{jk}$ are
Christoffel symbols. However for general Finsler metrics, $G^i$ are
not always quadratic in $y$. Thus  a natural  quantity is given by
\[ {\bf B}_{y} = B^i_{jkl}  \frac{\partial}{\partial x^i} \otimes dx^j  \otimes dx^k \otimes dx^l, \]
where
\[ B^i_{jkl} = B^i_{jkl}(x,y) := \frac{\partial^3 G^i}{\partial y^j \partial y^k \partial y^l}. \]
${\bf B}$ is called the {\it Berwald curvature}. Metrics with zero Berwald curvature are called {\it Berwald metrics}. It is obvious that all
Riemannian metrics are Berwald metrics.

  The {\it Landsberg curvature} ${\bf L}_{y} = L_{ijk} dx^i \otimes dx^j \otimes dx^k$  is defined by
\[ L_{ijk} = -\frac{1}{2} y^m g_{ml} B^l_{ijk} =-\frac{1}{2} y^m g_{ml}  \frac{\partial^3 G^l}{\partial y^i \partial y^j \partial y^k}.\]
Metrics with zero Landsberg curvature are called {\it Landsberg metrics}. It is obvious that all
Berwald metrics are Landsberg metrics. The {\it mean Landsberg curvature} ${\bf J}_y = J_k dx^k$ is defined by
\[ J_k = g^{ij} L_{ijk}. \]

 Related to weakly isotropic
  flag curvature we discussed here, S-curvature is first introduced by Z. Shen in \cite{Sh2} when he studied the
  volume comparison theorem in Finsler geometry. For a Finsler metric $F= F(x,y)$, the S-curvature is defined as the following.
  \[ {\bf S} = \frac{\partial G^m}{\partial y^m} - y^m \frac{\partial \ln \sigma_{F} }{\partial x^m}, \]
where $dV_{F} = \sigma_{F}(x) dx^1 ... dx^n$ is the Busemann-Hausdorff volume form. $F$ is said to have {\it almost isotropic
S-curvature} if there is a scalar function $c = c(x)$ such that
  \be \label{S-cur} {\bf S} = (n+1)(c F + \eta), \ee
where $\eta = \eta_i(x) y^i$ is a closed 1-form. $F$ is said to have {\it isotropic S-curvature} if $\eta =0$. Specially, $F$ is
said to have {\it constant S-curvature} if $\eta =0$ and $c$ is a constant. The following lemma shows the relationship between the
weakly isotropic flag curvature and almost S-curvature.
\begin{lem}[\cite{CMS}]\label{lemShen}
Let $F$ be a Finsler metric of scalar flag curvature ${\bf K} = {\bf K}(x,y)$. If $F$ has almost isotropic S-curvature as in
(\ref{S-cur}), then $F$ is of weakly isotropic flag curvature
\[ {\bf K} =\frac{3}{F} \frac{\partial c(x) }{\partial x^k} y^k +\sigma, \]
where $\sigma = \sigma(x)$ is a scalar function.
\end{lem}

 Projectively flat Finsler metrics compose  a large class of metrics of scalar flag curvature in Finsler geometry. A Finsler metric $F=F(x,y)$ on an open domain $\mathcal{U} \subset R^n$ is said to be \textit{projectively flat} in $\mathcal{U}$ if all geodesics are straight lines.
In 1903, G. Hamel proved that $F$ is projectively flat if and only if
\be \label{Hamel} \frac{\partial^2 F}{\partial x^l \partial y^k} - \frac{\partial^2 F}{\partial x^k \partial y^l} =0. \ee
This is also equivalent to its geodesic coefficients  $G^i = P(x,y) y^i$, where $P = F_{x^k}y^k/(2F) $ is called the {\it projective factor} of $F$. In this case, the flag curvature
 ${\bf K}$  is a scalar function on $T{\cal U}$ given by
\begin{equation*}
{\bf K} = \frac{P^2 - \frac{\partial P}{\partial x^m} y^m }{F^2}.\label{flagC}
\end{equation*}
In 1929, L. Berwald  proved the following lemma which plays an important role in our proof of Theorem \ref{thm2}.
\begin{lem}[\cite{Be2}]\label{lemBerwald}
Let $F=F(x,y)$ be a Finsler metric on an open subset $\mathcal{U} \subset R^n$. Then $F$ is projectively flat if and only if
there is a positively $y$-homogeneous function of degree one, $P = P(x, y)$, and a
positively homogeneous function of degree zero, $\mathbf{K} = \mathbf{K}(x, y)$, on $T \mathcal{U}  \simeq \mathcal{U} \times R^n$
such that
\begin{equation}\label{PF-F}
\frac{\partial F}{\partial x^k} = \frac{\partial (P F)}{\partial y^k},
\end{equation}
\begin{equation} \label{PF-P}
\frac{\partial P}{\partial x^k} = P \frac{\partial P}{\partial y^k} - \frac{1}{3 F} \frac{\partial (\mathbf{K} F^3)}{\partial y^k}.
\end{equation}
In this case, $P$ is the projective factor of $F$.
\end{lem}

\section{Berwald Connection and Main Lemmas} \label{sectionLemma}
In this section, we use Berwald connection to prove our main lemmas.
One can find details in \cite{ShenBook}. For simplicity, let
\[ F_{\cdot i} = \frac{\partial F}{\partial y^i}, \ \ \ F_{\cdot i\cdot j} = \frac{\partial^2 F}{\partial y^i\partial y^j}, \ \ \ F_{x^i} = \frac{\partial F}{\partial x^i}, \ \ \ .... \]
Let $"_|"$ denote the covariant derivative of Berwald connection.
The index "0" means the contraction with $y^i$. For example, for any scalar function $T=T(x,y)$ on $TM$
\[ T_{|0} = T_{|m}y^m, \ T_{|0\cdot k} = T_{|m\cdot k} y^m,\  T_{\cdot k|0} = T_{\cdot k|m}y^m,\  T_{|0|0} = T_{|l|m}y^l y^m,  \ \ \ .... \]
Let $R^{\ l}_{m\ ij}$ denote the hh-curvature of the Berwald connection. Then
\[ R^{\ l}_{m\ ij} = \frac{1}{3} \Big( \frac{\partial^2 R^l_{\ i}}{\partial y^m \partial y^j} - \frac{\partial^2 R^l_{\ j}}{\partial y^m \partial y^i}  \Big), \]
\[ R^{l}_{\ ij}  =y^m R^{\ l}_{m\ ij} = \frac{1}{3} \Big( \frac{\partial R^l_{\ i}}{ \partial y^j} - \frac{\partial R^l_{\ j}}{ \partial y^i}  \Big). \]

Using the Berwald connection,  the following Bianchi identities for
 $R^{\ l}_{m\ ij}$ are well-known.
\be \label{BianchiAdd}
 R^{\ l}_{i\ mk\cdot j} =  B^{l}_{ijk| m} - B^{l}_{ijm| k},
\ee

\be \label{Bianchi1}
R^{\ l}_{m\ ij|k} + R^{\ l}_{m\ jk|i} +R^{\ l}_{m\ ki|j} =
B^{l}_{mi u} R^{u}_{\ j k}+ B^{l}_{mj u} R^{u}_{\ k i}+ B^{l}_{mk u} R^{u}_{\ i j}.
\ee
Contracting the above equation with $y^m$ and $y^j$ yields
\be\label{Bianchi2}
R^{l}_{\ i|k} - R^{l}_{\ k|i} +R^{l}_{\ ki|0}=0.
\ee

Taking a trace of (\ref{Bianchi2}) over $l$ and $i$, we get
\be \label{Bianchi4}
R^{m}_{\ m|k} - R^{m}_{\ k|m} +R^{m}_{\ km|0} =0.
\ee
We  use the above two identities to prove the following two lemmas.
Now assume that Finsler metric $F$ is of scalar flag curvature ${\bf K} = {\bf K}(x,y)$. Then
\be\label{R^i_k}
R^i_{\ k} = {\bf K}F^2 h^i_k,
\ee
\be \label{R^m_m}
R^m_{\ m} =(n-1) {\bf K}F^2,
\ee
where \[ h^i_k = g^{ij}h_{jk}=\delta^i_k - \frac{1}{F} F_{\cdot k} y^i, \ \ \ h_{jk} = g_{jk} - F_{\cdot j} F_{\cdot k}. \]
By a direct computation, we have
\be \label{R^m_ij}
\begin{split}
R^m_{\ ij}=& \frac{1}{3}F^2({\bf K}_{\cdot j}h^m_i - {\bf K}_{\cdot i}h^m_j) - {\bf K}F(F_{\cdot i}\delta^m_j - F_{\cdot j}\delta^m_i) \\
=& \frac{F {\bf K}_{\cdot j}}{3}(F \delta^m_i- F_{\cdot i} y^m)-\frac{F {\bf K}_{\cdot i}}{3}(F \delta^m_j- F_{\cdot j} y^m)+ F {\bf K} (F_{\cdot j} \delta^m_i- F_{\cdot i} \delta^m_j),
\end{split}
\ee

\be \label{R^i_jkl}
\begin{split}
R^{\ i}_{j\ kl} =&  {\bf K}(g_{jl}\delta^i_k - g_{jk}\delta^i_l)+ \frac{F^2}{3}( {\bf K}_{\cdot j \cdot l}h^i_k -  {\bf K}_{\cdot j \cdot k}h^i_l ) + {\bf K}_{\cdot j}F( F_{\cdot l}\delta^i_k -  F_{\cdot k}\delta^i_l) \\
& + \frac{1}{3}{\bf K}_{\cdot l} (2 F F_{\cdot j}\delta^i_k - F F_{\cdot k} \delta^i_j - g_{jk}y^i)  - \frac{1}{3}{\bf K}_{\cdot k} (2 F F_{\cdot j}\delta^i_l - F F_{\cdot l} \delta^i_j - g_{jl}y^i),
\end{split}
\ee

\be\label{R^m_mij}
R^{\ m}_{m\ ij}=\frac{n+1}{3}F( {\bf K}_{\cdot i}F_{\cdot j} -{\bf K}_{\cdot j}F_{\cdot i}).
\ee

The following lemma was first proved in \cite{Najafi} when $F$ is of weakly isotropic flag curvature. Here we give the general version.
\begin{lem}\label{lemK2}
Let $F$ be a Finsler metric on an $n$-dimensional manifold $M$ ($n\geq3$) of scalar flag curvature ${\bf K}={\bf K}(x,y)$. Then
\be \label{K|k_1}
F {\bf K}_{|k} - F_{\cdot k} {\bf K}_{|0} - \frac{1}{3}F {\bf K}_{\cdot k|0} =0,
\ee
i.e.
\be\label{K|k_2}
\big[ \frac{{\bf K}_{|0}}{F} \big]_{\cdot k} = \frac{4}{3}\frac{{\bf K}_{\cdot k|0}}{F}.
\ee
\end{lem}
{\it Proof}: By a direct computation, we get
\[R^m_{\ m|k} =(n-1) {\bf K}_{|k}F^2, \ \ \ \  R^m_{\ k|m} = F^2 {\bf K}_{|k} - F F_{\cdot k}{\bf K}_{|0}, \]
\[ R^{m}_{\ km} = - \frac{n-2}{3}F^2{\bf K}_{\cdot k} -(n-1) {\bf K} F F_{\cdot k}, \]
\[ R^{m}_{\ km|j}y^j =- \frac{n-2}{3}F^2{\bf K}_{\cdot k|0} - (n-1)F F_{\cdot k}  {\bf K}_{|0}. \]
Substituting the above equations into (\ref{Bianchi4}) yields
\[ (n-2) \big\{ F {\bf K}_{|k} - F_{\cdot k} {\bf K}_{|0} - \frac{1}{3}F {\bf K}_{\cdot k|0} \big\} =0. \]
By the assumption $n\geq3$, we obtain (\ref{K|k_1}). Noting that  ${\bf K}_{|0\cdot k} = {\bf K}_{\cdot k|0}$, then
\[ [{\bf K}_{|0}  ]_{\cdot k} = {\bf K}_{|k} + {\bf K}_{\cdot k|0}. \]
Thus (\ref{K|k_1}) is equivalent to (\ref{K|k_2}).
\qed
\\

The following lemma plays an important role in proving our main results.
\begin{lem} \label{lemK1}
Let $F$ be a Finsler metric on $M$ of scalar flag curvature ${\bf K}={\bf K}(x,y)$. Then
\be \label{Kijk}
({\bf K}_{\cdot j |i}- {\bf K}_{\cdot i|j}) F_{\cdot k}+({\bf K}_{\cdot i |k} - {\bf K}_{\cdot k|i}) F_{\cdot j} + ( {\bf K}_{\cdot k |j}- {\bf K}_{\cdot j|k} ) F_{\cdot i} =0.
\ee
\end{lem}
{\bf Proof}: Differentiating (\ref{K|k_1}) with respect to $y^i$ yields
\be \label{K|k_i}
F_{\cdot i} {\bf K}_{|k} + F  {\bf K}_{|k\cdot i} - F_{\cdot k  \cdot i} {\bf K}_{|0} - F_{\cdot k} {\bf K}_{|i}  - F_{\cdot k } {\bf K}_{|0\cdot i}- \frac{1}{3}F_{\cdot i} {\bf K}_{\cdot k|0} - \frac{1}{3}F {\bf K}_{\cdot k|i}- \frac{1}{3}F {\bf K}_{\cdot k|0\cdot i} =0.
\ee

Exchanging the indices, we have
\be \label{K|i_k}
F_{\cdot k} {\bf K}_{|i} + F  {\bf K}_{|i\cdot k} - F_{\cdot k \cdot i} {\bf K}_{|0} - F_{\cdot i} {\bf K}_{|k}  - F_{\cdot i} {\bf K}_{|0\cdot k}- \frac{1}{3}F_{\cdot k} {\bf K}_{\cdot i|0} - \frac{1}{3}F {\bf K}_{\cdot i|k}- \frac{1}{3}F {\bf K}_{\cdot i|0\cdot k} =0.
\ee
Note that

\be
{\bf K}_{|k\cdot i} = {\bf K}_{\cdot i|k} - {\bf K}_{\cdot m} L_{lki} g^{lm}, \ \ {\bf K}_{|0 \cdot i} = {\bf K}_{\cdot i|0}
\ee
and
\be
{\bf K}_{\cdot i|0\cdot k} = {\bf K}_{\cdot i\cdot k|0} +  {\bf K}_{\cdot l} B^{l}_{ijk}y^j - {\bf K}_{\cdot i\cdot l} L_{mjk} g^{lm} y^j = {\bf K}_{\cdot i\cdot k|0} .
\ee

$\frac{1}{2 F }[(\ref{K|k_i}) - (\ref{K|i_k})] \times F_{\cdot j}$ yields
\[
( {\bf K}_{|i}- \frac{1}{3} {\bf K}_{|0\cdot i}) \frac{F_{\cdot j}F_{\cdot k}}{F} - ( {\bf K}_{|k} - \frac{1}{3}  {\bf K}_{|0\cdot k}  )\frac{F_{\cdot i}F_{\cdot j}}{F}+ \frac{2}{3} ( {\bf K}_{\cdot i|k} - {\bf K}_{\cdot k|i})F_{\cdot j} =0.
\]
The following two equations can be obtained similarly,
\[
( {\bf K}_{|j}- \frac{1}{3} {\bf K}_{|0\cdot j}) \frac{F_{\cdot k}F_{\cdot i}}{F} - ( {\bf K}_{|i} - \frac{1}{3}  {\bf K}_{|0\cdot i}  )\frac{F_{\cdot j}F_{\cdot k}}{F}+ \frac{2}{3} ( {\bf K}_{\cdot j|i} - {\bf K}_{\cdot i|j})F_{\cdot i} =0,
\]

\[
( {\bf K}_{|k}- \frac{1}{3} {\bf K}_{|0\cdot k}) \frac{F_{\cdot i}F_{\cdot j}}{F} - ( {\bf K}_{|j} - \frac{1}{3}  {\bf K}_{|0\cdot j}  )\frac{F_{\cdot k}F_{\cdot i}}{F}+ \frac{2}{3} ( {\bf K}_{\cdot k|j} - {\bf K}_{\cdot j|k})F_{\cdot i} =0.
\]
By summing up the above three equations, we obtain (\ref{Kijk}).
\qed
\\

Based on  Lemma \ref{lemK1}, when ${\bf K} = \frac{3\theta}{F} +\sigma$, we can prove that $\theta$ must be closed.
\begin{lem}\label{lem1}
Let $F$ be a Finsler metric on an $n$-dimensional manifold $M$ ($n\geq 3$) of weakly isotropic flag curvature
\[ {\bf K} = \frac{3\theta}{F} + \sigma,\]
where $\theta = \theta_i(x) y^i$ is a $1$-form and $\sigma=\sigma(x)$ is a scalar function. Then
$\theta$ must be a closed $1$-form.
\end{lem}
{\it Proof}:
By a direct computation, we get
\[K_{\cdot i} = \frac{3\theta_i}{F} - \frac{3\theta}{F^2} F_{\cdot i},\]
\be \label{K_ij} K_{\cdot i|j} = \frac{3\theta_{i|j}}{F} - \frac{3\theta_{|j}}{F^2} F_{\cdot i}.  \ee
Plugging (\ref{K_ij}) into (\ref{Kijk}) yields

\be \label{FiFjFk}
 ( \theta_{k |j} - \theta_{j |k} )F_{\cdot i} + (\theta_{i |k} - \theta_{k |i})F_{\cdot j} + (\theta_{j |i}- \theta_{i |j} )F_{\cdot k} =0.
\ee
By $\theta_{i|j} = \theta_{ix^j} - \theta_m G^m_{\cdot i\cdot j}$, $\theta_{j |i}- \theta_{i |j} =\theta_{j x^i}- \theta_{i x^j}$ and $ F_{\cdot i} = \Big[ \frac{1}{2}F^2 \Big]_{\cdot i} $, we have
\be \label{FiFjFk}
 ( \theta_{k x^j} - \theta_{j x^k} )\Big[ \frac{1}{2}F^2\Big]_{\cdot i} + (\theta_{i x^k} - \theta_{k x^i})\Big[ \frac{1}{2}F^2\Big]_{\cdot j} + (\theta_{j x^i}- \theta_{i x^j} )\Big[\frac{1}{2}F^2\Big]_{\cdot k} =0.
\ee
Noticing $g_{ij} = \frac{1}{2}[F^2]_{\cdot i \cdot j}$
 and differentiating (\ref{FiFjFk}) respect to $y^l$ yields
\be \label{gilgjlgkl}
 ( \theta_{k x^j} - \theta_{j x^k} )g_{il} + (\theta_{i x^k} - \theta_{k x^i})g_{jl} + (\theta_{j x^i}- \theta_{i x^j} )g_{kl} =0.
\ee
Contracting the above equation with $g^{kl}$, we get
\be
(n-2) (\theta_{j x^i}- \theta_{i x^j} )=0.
\ee
By the assumption $ n\geq 3 $, then $\theta_{j x^i}= \theta_{i x^j}$. Thus $\theta$ is a closed $1$-form.
\qed

By Lemma \ref{lemK1} and \ref{lemK2}, we can prove the following main lemma which  plays an important role in our proof of Theorem \ref{thm1}.
\begin{lem} \label{mainlem}
Let $F$ be a Finsler metric on an $n$-dimensional manifold $M$ ($n\geq 3$) of weakly isotropic flag curvature
\[ {\bf K} = \frac{3\theta}{F} + \sigma,\]
where $\theta = \theta_i(x) y^i$ is a $1$-form and $\sigma=\sigma(x)$ is a scalar function. Then
$\theta$ is a closed 1-form and there exists a scalar function $f = f(x)$ such that
\be \label{mainlemEq}
f F^2 - \sigma_{|0} F - \theta_{|0} =0.
\ee
In this case,
\be \label{mainlemEq_i}
\theta_{|i} = f F F_{\cdot i} -\frac{1}{2} \sigma_{|i} F - \frac{1}{2}\sigma_{|0} F_{\cdot i},
\ee
\be \label{mainlemEq_ivj}
\theta_{|i|j} = f_{|j} F F_{\cdot i} -\frac{1}{2} \sigma_{|i|j} F - \frac{1}{2}\sigma_{|0|j} F_{\cdot i}.
\ee

\end{lem}
{\it Proof}: By Lemma \ref{lem1}, $\theta$ is closed. Then
\be \label{thetaclosed1}
 \theta_{m|k} - \theta_{k|m} =0, \ \ \  \theta_{|k} - \theta_{k|m}y^m =0, \ee
and
\be  \label{thetaclosed2}
( \theta_{|m}y^m)_{\cdot k} = 2 \theta_{|k}. \ee
By a direct computation, we have
\[ {\bf K}_{|0} = \frac{3\theta_{|0}}{F} + \sigma_{|0}, \]
\[ {\bf K}_{\cdot k|0} =  \frac{3\theta_{\cdot k|0}}{F} - \frac{3 \theta_{|0}}{F^2}F_{\cdot k} = \frac{3\theta_{|k}}{F} - \frac{3 \theta_{|0}}{F^2}F_{\cdot k} = \frac{3}{2} F \big[ \frac{\theta_{|0}}{F^2} \big]_{\cdot k}. \]
Plugging the above two equations into (\ref {K|k_2}) yields

\be
\big[  \frac{\sigma_{|0} }{F} + \frac{\theta_{|0} }{F^2}  \big]_{\cdot k} =0.
\ee
Then there exists a scalar function $f = f(x)$ such that (\ref{mainlemEq}) holds.
Differentiating (\ref{mainlemEq}) respect to $y^i$ yields
\[ 2f F F_{\cdot i} - \sigma_{|i} F - \sigma_{|0} F_{\cdot i} - 2 \theta_{|i} =0. \]
Differentiating the above equation along the direction $\frac{\delta }{\delta x^j}$ yields
\[ 2f_{|j} F F_{\cdot i} - \sigma_{|i|j} F - \sigma_{|0|j} F_{\cdot i} - 2 \theta_{|i|j} =0. \]
Then we obtain (\ref{mainlemEq_i}) and (\ref{mainlemEq_ivj}).
\qed

Our proofs of the Theorem \ref{thm1} and \ref{thm2} are based on the following lemma.
\begin{lem}\label{lemtrivial}
Let $F = F(x,y)$ be a Finsler metric. If there is a scalar function $a = a(x)$, a 1-form $\eta = \eta_i(x)y^i$ and a 2-form $\xi=\xi_{ij} y^i y^j$ such that
\be \label{lemtrivialeq} a F^2 + \eta F + \xi =0. \ee
If $\xi \neq 0$ or $\eta \neq 0$,
then $F$ is a Randers metric.
\end{lem}
{\it Proof}: Obviously, (\ref{lemtrivialeq})  is a quadratic equation of $F$. By the assumption that  $\xi \neq 0$ or $\eta \neq 0$, then
\[ a \neq 0 \]
and
\[ F = \frac{-\eta + \sqrt{\eta^2 - 4 a \xi}}{2a} \]
is a Randers metric.
\qed

\section{Proof of Theorem \ref{thm1}}
In this section, we prove Theorem \ref{thm1} by Proposition \ref{prop1} and Lemma \ref{lemlast}.
By Lemma \ref{mainlem}, the following proposition can be proved now. The key idea is to use the Ricci identity several times.
\begin{prop}\label{prop1}
Let $F$ be a Finsler metric on an $n$-dimensional manifold $M$ ($n\geq 3$) of weakly isotropic flag curvature
\[ {\bf K} = \frac{3\theta}{F} + \sigma,\]
where $\theta = \theta_i(x) y^i \neq 0 $ is a $1$-form and $\sigma=\sigma(x)$ is a scalar function. Then there exists a scalar function
$\lambda =\lambda(x)$ such that $\theta = \lambda \sigma_{|0}$ and  one of the following holds.

(i) $F$ is a Randers metric;

(ii) The scalar function $\sigma$ and $\lambda$ satisfy
\be  \label{thmeq_last1} 4\lambda(f_{|0} +\sigma \theta) + \sigma_{|0} =0 \ee
and
\be  \label{thmeq_last2} 6\lambda^3 \sigma_{|0} + \lambda_{|0}  =0, \ee
where $f=f(x)$ satisfies (\ref{mainlemEq}).
\end{prop}
{\it Proof}:
 By the assumption ${\bf K} = \frac{3\theta}{F} + \sigma$, plugging it into (\ref{R^m_ij}) yields
\[ R^m_{\ ij} = \frac{1}{F} (\theta_j \delta^m_i - \theta_i \delta^m_j) - y^m (\theta_j F_{\cdot i}- \theta_i F_{\cdot j}) -\frac{\theta}{F} (F_{\cdot j} \delta^m_i - F_{\cdot i}\delta^m_j) + F {\bf K} (F_{\cdot j}\delta^m_i - F_{\cdot i}\delta^m_j). \]
By the Ricci identity we have
\[
\theta_{|i|j} - \theta_{|j|i} = \theta_m R^m_{\ ij} = F {\bf K} (\theta_i F_{\cdot j} - \theta_j F_{\cdot i}).
  \]
Substituting (\ref{mainlemEq_ivj}) into the above equation, then by $\sigma_{|i|j} = \sigma_{|j|i}$ we get
\be \label{fF1}
f_{|j} F F_{\cdot i} - f_{|i} F F_{\cdot j} -\frac{1}{2}\sigma_{|0|j} F_{\cdot i} + \frac{1}{2}\sigma_{|0|i} F_{\cdot j} =  F {\bf K} (\theta_i F_{\cdot j} - \theta_j F_{\cdot i}).
\ee
Contracting (\ref{fF1}) with $y^j$ yields
\be \label{fF2}
f_{|0} F F_{\cdot i} - f_{|i} F^2 -\frac{1}{2}\sigma_{|0|0} F_{\cdot i} + \frac{1}{2}\sigma_{|0|i} F =  F {\bf K} (\theta_i F - \theta F_{\cdot i}).
\ee
Noting that
\[ f_{|0} F F_{\cdot i} - f_{|i} F^2 = - F^3 \big[ \frac{f_{|0}}{F}\big]_{\cdot i}, \ \ \
\sigma_{|0|i} F - \sigma_{|0|0} F_{\cdot i} = \frac{1}{2} F^3 \big[ \frac{\sigma_{|0|0}}{F^2}\big]_{\cdot i}
\]
and
\[ (\theta_i F - \theta F_{\cdot i}) = F^2 \big[ \frac{\theta}{F} \big]_{\cdot i} =  \frac{1}{3}F^2 \big[ \frac{3 \theta}{F} + \sigma\big]_{\cdot i} = \frac{1}{3}F^2 {\bf K}_{\cdot i},  \]
we can rewrite (\ref{fF2})  into
\[
F^3 \big[ \frac{f_{|0}}{F} -  \frac{1}{4} \frac{\sigma_{|0|0}}{F^2} + \frac{1}{6} {\bf K}^2 \big]_{\cdot i} = F^3 \big[ \frac{f_{|0}}{F} + \frac{ \sigma\theta}{F} -  \frac{1}{4} \frac{\sigma_{|0|0}}{F^2} + \frac{3\theta^2}{2F^2}  \big]_{\cdot i}=0.
\]
Then there exists a scalar function $h = h(x)$ such that
\[\frac{f_{|0}}{F} + \frac{ \sigma\theta}{F} -  \frac{1}{4} \frac{\sigma_{|0|0}}{F^2} + \frac{3\theta^2}{2F^2}   = - \frac{1}{2}h. \]
Which is equivalent to
\be \label{sigma_v0v0}
\frac{1}{2} \sigma_{|0|0} = h  F^2 + 2(f_{|0} +\sigma \theta)F + 3 \theta^2.
\ee
Now we divided our left proof into the following two cases.

{\bf Case (a) $\sigma_{|0} =0$.} In this case, (\ref{sigma_v0v0}) becomes into
\[ h F^2 + 2(f_{|0} +\sigma \theta)F + 3 \theta^2 = 0. \]
It is easy to see that $F$ must be a Randers metric by Lemma \ref{lemtrivial}.

{\bf Case (b) $\sigma_{|0} \neq 0$.} In this case, we prove that $\theta$ is parallel to $\sigma_{|0}$.
Differentiating (\ref{sigma_v0v0}) respect to $y^i$ yields
\be \label{sigma_v0vi}
 \sigma_{|0|i} = 2 h F F_{\cdot i} + 2(f_{|i} +\sigma \theta_i)F + 2(f_{|0} +\sigma \theta)F_{\cdot i} + 6 \theta \theta_i.
\ee
Differentiating (\ref{sigma_v0vi}) along the direction $\frac{\delta }{\delta x^j}$ yields
\[
 \sigma_{|0|i|j} = 2h_{|j} F F_{\cdot i} + 2(f_{|i|j} +\sigma_{|j} \theta_i + \sigma \theta_{i|j} )F + 2(f_{|0} +\sigma \theta)_{|j}F_{\cdot i} + 6 \theta_{|j} \theta_i + 6\theta \theta_{i|j}.
\]
Then
\be \label{sigma_vivj-vjvi}
\begin{split}
& \sigma_{|0|i|j} - \sigma_{|0|j|i} \\
= & 2 ( h_{|j} F_{\cdot i} - h_{|i} F_{\cdot j})F + 2 (\sigma_{|j} \theta_i -\sigma_{|i} \theta_j)F \\
&+ 2(f_{|0} +\sigma \theta)_{|j}F_{\cdot i} - 2(f_{|0} +\sigma \theta)_{|i}F_{\cdot j} +  6 (\theta_{|j} \theta_i -\theta_{|i} \theta_j).
\end{split}
\ee
Here we used $f_{|i|j} = f_{|i|j}$ and $\theta_{i|j} = \theta_{j|i}$.
Contracting (\ref{sigma_vivj-vjvi}) with $y^j$ yields
\be \label{sigma_vivj-vjviyj}
\begin{split}
& \sigma_{|0|i|0} - \sigma_{|0|0|i} \\
= & 2 ( h_{|0} F_{\cdot i} - h_{|i} F )F + 2 (\sigma_{|0} \theta_i -\sigma_{|i} \theta)F \\
&+ 2(f_{|0} +\sigma \theta)_{|0} F_{\cdot i} - 2(f_{|0} +\sigma \theta)_{|i}F  +  6 (\theta_{|0} \theta_i -\theta_{|i} \theta) \\
= & - F^3 \big[  \frac{2 h_{|0}}{F} + \frac{(f_{|0} +\sigma \theta)_{|0}}{F^2} \big]_{\cdot i}  +  2 (\sigma_{|0} F \theta_i -\sigma_{|i} F \theta) +  6 (\theta_{|0} \theta_i -\theta_{|i} \theta).
\end{split}
\ee
By (\ref{mainlemEq}) and (\ref{mainlemEq_i}), we have
\[ \sigma_{|0} F = f F^2 - \theta_{|0}, \]
\[  \sigma_{|i} F = 2 f F F_{\cdot i} - \sigma_{|0} F_{\cdot i} - 2 \theta_{|i} = f F F_{\cdot i} + \frac{1}{F}\theta_{|0} F_{\cdot i} - 2 \theta_{|i}.\]
Substituting the above two equations into (\ref{sigma_vivj-vjviyj}) yields
\be \label{sigma_vivj-vjviyj_1}
\begin{split}
 & \sigma_{|0|i|0} - \sigma_{|0|0|i}\\
= & - F^3 \big[  \frac{2 h_{|0}}{F} + \frac{(f_{|0} +\sigma \theta)_{|0}}{F^2} \big]_{\cdot i} + 2(f F^2\theta_i - f F F_{\cdot i}\theta ) - \frac{2}{F}\theta_{|0} F_{\cdot i}\theta + 4 \theta_{|0} \theta_i - 2\theta_{|i} \theta \\
= & - F^3 \big[  \frac{2 h_{|0}}{F} + \frac{(f_{|0} +\sigma \theta)_{|0}}{F^2} - 2 f \frac{\theta}{F}\big]_{\cdot i} - \frac{2}{F}\theta_{|0} F_{\cdot i} \theta + 4 \theta_{|0} \theta_i - 2\theta_{|i} \theta.
\end{split}
\ee
On the other hand, by Ricci identity we have
\be \label{sigma_vivj-vjviyj_2}
\begin{split}
\sigma_{|0|i|0} - \sigma_{|0|0|i}  = \sigma_{|m}R^m_{\ ij}y^j =& 2 (\frac{3\theta}{F}+ \sigma)( \theta_{|0} F_{\cdot i} - F\theta_{|i})\\
 =& \frac{6}{F}  \theta_{|0} F_{\cdot i} \theta - 6 \theta_{|i}\theta -  F^3 \big[ \sigma \frac{\theta_{|0} }{F^2} \big]_{\cdot i}.
 \end{split}
\ee

By (\ref{sigma_vivj-vjviyj_1}) and (\ref{sigma_vivj-vjviyj_2}), we have
\be
\begin{split}
& \big[  \frac{2h_{|0}}{F} + \frac{(f_{|0} +\sigma \theta)_{|0}}{F^2} - 2 f \frac{\theta}{F} - \sigma \frac{\theta_{|0}}{F^2}\big]_{\cdot i}\\
= &  \frac{4}{F^3} (  \theta_{|0} \theta_i + \theta_{|i} \theta - \frac{2}{F}\theta_{|0} F_{\cdot i} \theta)
=  4 \frac{\theta_{|0}}{F^2} \big[ \frac{\theta}{F} \big]_{\cdot i} + 2 \big[ \frac{\theta_{|0}}{F^2} \big]_{\cdot i}  \frac{\theta}{F}
\end{split}
\ee
Then
\be \label{thmeqstar}
\big[  \frac{2 h_{|0}}{F} + \frac{(f_{|0} +\sigma \theta)_{|0}}{F^2} - 2 f \frac{\theta}{F} - \sigma \frac{\theta_{|0}}{F^2} - 2 \frac{\theta_{|0}}{F^2}  \frac{\theta}{F} \big]_{\cdot i} = 2 \frac{\theta_{|0}}{F^2} \big[ \frac{\theta}{F} \big]_{\cdot i}.
\ee
Thus, we obtain
\[
\big[ \frac{\theta_{|0}}{F^2}\big]_{\cdot j} \big[ \frac{\theta}{F} \big]_{\cdot i} = \big[ \frac{\theta_{|0}}{F^2}\big]_{\cdot i} \big[ \frac{\theta}{F} \big]_{\cdot j}.
\]
By (\ref{mainlemEq}), the above equation  is equivalent to
\[
\big[ \frac{\sigma_{|0}}{F}\big]_{\cdot j} \big[ \frac{\theta}{F} \big]_{\cdot i} = \big[ \frac{\sigma_{|0}}{F}\big]_{\cdot i} \big[ \frac{\theta}{F} \big]_{\cdot j}.
\]
Simplifying the above equation yields
\be \label{thetasigma}
(\sigma_{|i} \theta_j -\sigma_{|j} \theta_i )F^2 + \sigma_{|m}y^m (\theta_i F F_{\cdot j} - \theta_j F F_{\cdot i}) + \theta (\sigma_{|j}FF_{\cdot i} - \sigma_{|i} F F_{\cdot j}) =0.
\ee
Differentiating (\ref{thetasigma}) with respect to $y^k$ and $y^l$ yields
\be \label{thetasigma1}
\begin{split}
0= & 2 (\sigma_{|i} \theta_j -\sigma_{|j} \theta_i )g_{kl} + \sigma_{|k} (\theta_i g_{jl} - \theta_j g_{il})
+ \sigma_{|l} (\theta_i g_{jk} - \theta_j g_{ik}) \\
&+ 2 \sigma_{|m}y^m (\theta_i C_{jkl} - \theta_j C_{ikl})
+ \theta_{k} (\sigma_{|j}g_{il} - \sigma_{|i} g_{jl}) \\
&+ \theta_{l} (\sigma_{|j}g_{ik} - \sigma_{|i} g_{jk})
+ 2 \theta (\sigma_{|j}C_{ikl} - \sigma_{|i}C_{jkl}),
\end{split}
\ee
where $C_{ijk} = \frac{1}{4}[F^2]_{\cdot i \cdot j \cdot k}$ and $I_i = g^{jk}C_{ijk}$.
Contracting (\ref{thetasigma1}) with $g^{kl}$, we get
\be \label{thetasigma2}
\begin{split}
0= ( n -2) (\sigma_{|i} \theta_j -\sigma_{|j} \theta_i )+  \sigma_{|0} (\theta_i I_{j} - \theta_j I_{i})
+  \theta (\sigma_{|j}I_{i} - \sigma_{|i}I_{j}).
\end{split}
\ee
By $I_j y^j =0$, contracting (\ref{thetasigma2}) with $y^j$ yields
\[ 0= (n -2) (\sigma_{|i} \theta -\sigma_{|0} \theta_i ). \]
Thus, there exists a scalar function $\lambda=\lambda(x)$ such that
\be\label{thetalambdasigma} \theta_i =  \lambda \sigma_{|i}. \ee

By a direct computation, we get
\be \label{thmcase2}
\begin{split}
 \theta_{|0} =& \lambda_{|0} \sigma_{|0} + \lambda \sigma_{|0|0}
                =  2 \lambda h F^2 + 4 \lambda (f_{|0} +\sigma \theta)F + 6 \lambda \theta^2 + \lambda_{|0} \sigma_{|0}.
 \end{split}
\ee
Here the second equality is from (\ref{sigma_v0v0}). Substituting (\ref{thmcase2}) back into (\ref{mainlemEq}) yields
\be \label{thmeqfh}
 (f -2\lambda h) F^2 - [4\lambda(f_{|0} +\sigma \theta) + \sigma_{|0}]F - ( 6\lambda^3 \sigma_{|0} + \lambda_{|0}) \sigma_{|0} =0.
 \ee
If $ 4\lambda(f_{|0} +\sigma \theta) + \sigma_{|0} \neq 0$ or $6\lambda^3 \sigma_{|0} + \lambda_{|0}  \neq 0$, then by Lemma \ref{lemtrivial}, $F$ is a Randers metric. Otherwise (\ref{thmeq_last1}) and (\ref{thmeq_last2}) hold.
\qed

In the following lemma we prove that the case (ii) in Proposition
\ref{prop1} is impossible if the metric is not a Randers metric.
\begin{lem} \label{lemlast}
Let $F$ be a  Finsler metric on an $n$-dimensional manifold $M$ ($n\geq 2$) of weakly isotropic flag curvature
\[ {\bf K} = \frac{3\theta}{F} + \sigma,\]
where $\theta = \theta_i(x) y^i = \lambda \sigma_{|i} y^i $ is a closed $1$-form and $\sigma=\sigma(x)$, $\lambda =\lambda(x)$ are two scalar functions.  Assume that $F$ is not a Randers metric.
If
\be \label{lemlast1} 4\lambda(f_{|0} +\sigma \theta) + \sigma_{|0} =0 \ee
and
\be \label{lemlast2} 6\lambda^3 \sigma_{|0} + \lambda_{|0}  =0, \ee
where $f=f(x)$ satisfies (\ref{mainlemEq}), then
\[ \theta = \lambda \sigma_{|i} y^i =0. \]
\end{lem}
{\it Proof}: If at some point $x_o$, $\lambda(x_o)=0$, then at this point $\theta=0$. Thus we only need to consider the point $x$ such that $\lambda(x)\neq 0$.
Plugging  $\theta = \lambda \sigma_{|i} y^i$ back to (\ref{mainlemEq_i}) yields
\[ \theta_{|0} = f F^2 - \frac{\theta F}{\lambda}. \]
Then by a direct substitution, (\ref{thmeqstar}) can be written into
\be \label{thmeqstar_new}
\big[  \frac{2 h_{|0}}{F} + \frac{(f_{|0} +\sigma \theta)_{|0}}{F^2} - 4 f \frac{\theta}{F} - \sigma \frac{\theta_{|0}}{F^2} - 2 \frac{\theta_{|0}}{F^2}  \frac{\theta}{F}  - \frac{\theta^2}{\lambda F^2} \big]_{\cdot i} =0.
\ee
This means that there is a function $r=r(x)$ such that
\be \label{r}
\frac{2 h_{|0}}{F} + \frac{(f_{|0} +\sigma \theta)_{|0}}{F^2} - 4 f \frac{\theta}{F} - \sigma \frac{\theta_{|0}}{F^2} - 2 \frac{\theta_{|0}}{F^2}  \frac{\theta}{F}  - \frac{\theta^2}{\lambda F^2} = r
\ee

By (\ref{lemlast1}), we have
\be
\begin{split}
(f_{|0} +\sigma \theta)_{|0} & = - (\frac{\sigma_{|0}}{4\lambda})_0 =  -\frac{1}{4}( \frac{\sigma_{|0|0}}{\lambda}- \frac{ \sigma_{|0} \lambda_{|0}}{\lambda^2}) \\
& = -\frac{1}{16\lambda^2}(F-4 \lambda_{|0} )\sigma_{|0} - \frac{1}{2\lambda} (h F^2 + 3\theta^2).
\end{split}
\ee
The last equality is from (\ref{sigma_v0v0}) and (\ref{lemlast1}). Then by the above expression of $(f_{|0} +\sigma \theta)_{|0}$ and  (\ref{lemlast2}), (\ref{r}) can be simplified into
\[
-\frac{1}{2} \lambda (2 f \lambda \sigma + 2r
\lambda + h) F^2+(-6 f \lambda^3 \sigma_{|0} + \lambda^2 \sigma \sigma_{|0}+ 2 h_{|0} \lambda^2- \frac{1}{16} \sigma_{|0}) F+\lambda^3 (\lambda^2-1)\sigma_{|0}^2 =0.
\]
If $F$ is not a Randers metric, then by Lemma \ref{lemtrivial}, $ \lambda^3 (\lambda^2-1)\sigma_{|0}^2 \neq 0$.
Then $\lambda =const.$ or $\sigma_{|0}=0$. In fact, if $\lambda$ is a constant, then by  (\ref{lemlast2}), $\sigma_{|0}=0$.
Thus $\sigma_{|0} =0$.\qed

{\bf Proof of Theorem \ref{thm1}}: By the assumption, $F$ is of weakly isotropic flag curvature
\[ {\bf K} = \frac{3\theta}{F} + \sigma.\]
If $\theta=0$, then by Schur's Lemma $F$ is of constant flag curvature. If $\theta \neq 0$, then by Proposition \ref{prop1} and Lemma \ref{lemlast} $F$ must be a Randers metric. Further, by the results in \cite{ChSh} $F$ must be expressed in (\ref{Randers}).
\qed
\section{Projectively flat Finsler metrics}
Without the restriction on the dimension of the manifold,  we study projectively flat Finsler metrics of weakly isotropic flag curvature in this section. The following lemma shows the equation which the flag curvature should satisfy.

\begin{lem}
Let $F$ be a projectively flat Finsler metric of scalar flag curvature ${\bf K} = {\bf K}(x,y)$. Then
\be\label{ProjK}
\frac{F}{3}( {\bf K}_{\cdot l}P_{\cdot k}-{\bf K}_{\cdot k}P_{\cdot l}) + P ({\bf K}_{\cdot l}F_{\cdot k}-{\bf K}_{\cdot k}F_{\cdot l})
 + {\bf K}_{x^k}F_{\cdot l}-{\bf K}_{x^l}F_{\cdot k} + \frac{F}{3}({\bf K}_{\cdot l x^k}-{\bf K}_{\cdot k x^l})=0,
\ee
where $P=P(x,y)$ is the projective factor of $F$.
\end{lem}
{\it Proof}:
By (\ref{PF-P}) and a direct computation we have
\be \label{PF-Pkxl}
P_{ x^l \cdot k} = \frac{1}{2}[P^2]_{\cdot k\cdot l} - F F_{\cdot l} {\bf K}_{\cdot k} - \frac{1}{2}{\bf K} [F^2]_{\cdot k\cdot l} - \frac{2}{3} F F_{\cdot k} {\bf K}_{\cdot l} -\frac{1}{3}F^2 {\bf K}_{\cdot k\cdot l}
\ee

\be
P_{x^k x^l} = P_{x^l} P_{\cdot k} + P P_{k x^l}- F{\bf K}_{x^l} F_{\cdot k}- {\bf K}F_{x^l}F_{\cdot k} - {\bf K} F F_{\cdot k x^l} - \frac{2}{3}F F_{x^l} {\bf K}_{\cdot k} -\frac{1}{3}F^2 {\bf K}_{\cdot k x^l}.
\ee
Plugging (\ref{PF-F}), (\ref{PF-P}) and (\ref{PF-Pkxl}) into above equation yields
\be
\begin{split}
P_{x^k x^l} =& 2 P P_{\cdot k} P_{\cdot l} + P^2 P_{\cdot k \cdot l}- {\bf K} P F_{\cdot k}F_{\cdot l} - {\bf K} F (P_{\cdot k} F_{\cdot l}+ P_{\cdot l}F_{\cdot k}) - \frac{{\bf K} P}{2}[F^2]_{\cdot k\cdot l} \\
& - \frac{P F^2}{3}{\bf K}_{\cdot k\cdot l} - {\bf K} F F_{\cdot k x^l} - \frac{F^2}{3}({\bf K}_{\cdot k}P_{\cdot l}+{\bf K}_{\cdot l}P_{\cdot k}) - \frac{2 P F}{3}({\bf K}_{\cdot k}F_{\cdot l} + {\bf K}_{\cdot l}F_{\cdot k})\\
& -\frac{F^2}{3}{\bf K}_{\cdot k}P_{\cdot l} - P F{\bf K}_{\cdot k}F_{\cdot l} - F {\bf K}_{x^l}F_{\cdot k} - \frac{F^2}{3}{\bf K}_{\cdot k x^l}.
\end{split}
\ee
Then
\be
\begin{split}
0 = P_{x^k x^l} - P_{x^l x^k}
=&\frac{F^2}{3}( {\bf K}_{\cdot l}P_{\cdot k}-{\bf K}_{\cdot k}P_{\cdot l}) + P F({\bf K}_{\cdot l}F_{\cdot k}-{\bf K}_{\cdot k}F_{\cdot l})\\
 &+ F({\bf K}_{x^k}F_{\cdot l}-{\bf K}_{x^l}F_{\cdot k}) + \frac{F^2}{3}({\bf K}_{\cdot l x^k}-{\bf K}_{\cdot k x^l}).
\end{split}
\ee
Here G. Hamel's equation (\ref{Hamel}) is used.
\qed

Now we can prove that when
\[ {\bf K} = \frac{3\theta}{F} + \sigma, \]
$\theta$ must be closed. By the equations (\ref{PF-F}) and  (\ref{PF-P}) in Lemma \ref{lemBerwald}, we get the following equations.
\be \label{ProjK1}
{\bf K}_{\cdot k} = \frac{3}{F^2}( \theta_k F - \theta F_{\cdot k} ),
\ee
\be \label{ProjK2}
{\bf K}_{x^k} = \frac{3}{F^2}(\theta_{x^k} F - \theta F_{x^k})+\sigma_{x^k} = \frac{3}{F^2}\Big[ \theta_{x^k} F - \theta (P_{\cdot k}F+ P F_{\cdot k})\Big]+\sigma_{x^k},
\ee
\be \label{ProjK3}
{\bf K}_{x^l}y^l = \frac{3}{F}(\theta_{x^l}y^l  - 2 \theta P )+\sigma_{x^l}y^l,
\ee
\be \label{ProjK4}
{\bf K}_{\cdot k x^l}y^l = \frac{12 P}{F^2}(\theta F_{\cdot k} - \theta_k F) + \frac{3}{F^2} (\theta_{k x^l}y^l F + 2 PF \theta_k - \theta_{x^l}y^l F_{\cdot k} - \theta F_{x^k}).
\ee
Then we can prove the following lemma.
\begin{lem} \label{lemthetaclose}
Let $F$ be a projectively flat Finsler metric of weakly isotropic flag curvature
\[ {\bf K} = \frac{3\theta}{F} + \sigma,\]
where $\theta = \theta_i(x) y^i$ is a $1$-form and $\sigma=\sigma(x)$ is a scalar function. Then
$\theta$ must be a closed $1$-form. In this case, there exists a scalar function $a= a(x)$ such that
\be \label{aF}
a F^2 - \sigma_{x^l}y^l F + 2\theta P - \theta_{x^l}y^l =0.
\ee
\end{lem}
{\it Proof}: Contracting (\ref{ProjK}) with $y^l$ yields
\be
- \frac{4 P F}{3}{\bf K}_{\cdot k}
 +F {\bf K}_{x^k} -{\bf K}_{x^l}y^l F_{\cdot k} - \frac{F}{3}{\bf K}_{\cdot k x^l}y^l=0.
\ee
Substituting (\ref{ProjK1})-(\ref{ProjK4}) into the above equation yields
\[
\frac{2}{F}(2 \theta P F_{\cdot k} - \theta F P_{\cdot k} - P F \theta_k) + F \sigma_{x^k} -\sigma_{x^l} y^l F_{\cdot k} - \frac{2}{F} \theta_{x^l}y^l F_{\cdot k}-\theta_{k x^l}y^l + 3 \theta_{x^k} =0.
\]
The above equation can be written into
\be \label{diffEQ}
2 \frac{\theta_{k x^l}y^l - \theta_{x^k}}{F^2} = \Big[ -  \frac{2 \theta P}{F^2} + \frac{\sigma_{x^l}y^l}{F} + \frac{\theta_{x^l}y^l} {F^2}  \Big]_{\cdot k}.
\ee
Then we have
\be
\Big[ \frac{\theta_{k x^l}y^l - \theta_{x^k}}{F^2} \Big]_{\cdot j} = \Big[ \frac{\theta_{j x^l}y^l - \theta_{x^j}}{F^2} \Big]_{\cdot k}.
\ee
By a direct computation, we get
\[
(\theta_{k x^j} - \theta_{j x^k}) F^2 - (\theta_{k x^l}y^l - \theta_{x^k})\Big[\frac{1}{2}F^2 \Big]_{\cdot j} + (\theta_{j x^l}y^l -\theta_{x^j})\Big[ \frac{1}{2}F^2 \Big]_{\cdot k} =0.
\]
Differentiating the above equation respect to $y^i$ yields
\[
2 (\theta_{k x^j} - \theta_{j x^k}) F F_{\cdot i} - (\theta_{k x^i} - \theta_{i x^k})F F_{\cdot j} + (\theta_{j x^i} -\theta_{i x^j})F F_{\cdot k}- (\theta_{k x^l}y^l - \theta_{x^k})g_{ij} + (\theta_{j x^l}y^l -\theta_{x^j})g_{ki} =0.
\]
Contracting it with $g^{ij}$, we get
\[
n(\theta_{x^k} - \theta_{kx^l}y^l)=0.
\]
Then
\[ \theta_{l x^k} = \theta_{k x^l}. \]
Thus $\theta$ is a closed 1-form. Substituting it back into (\ref{diffEQ}) yields
\[ \Big[ -  \frac{2 \theta P}{F^2} + \frac{\sigma_{x^l}y^l}{F} + \frac{\theta_{x^l}y^l} {F^2}  \Big]_{\cdot k} =0. \]
Then we obtain (\ref{aF}).
\qed

\begin{lem} \label{lemprobF}
Let $F$ be a projectively flat Finsler metric of weakly isotropic flag curvature
\[ {\bf K} = \frac{3\theta}{F} + \sigma,\]
where $\theta = \theta_i(x) y^i$ is a $1$-form and $\sigma=\sigma(x)$ is a scalar function. Then one of the following holds. \\
(i)  $F$ is a Randers metric; \\
(ii) $\sigma = \sigma(x) \neq const.$ and  there exist two scalar function $b = b(x)$ and $\lambda =\lambda (x)$ such that
\be  \label{lembF}
 \frac{1}{2}\sigma_{x^k x^l}y^k y^l - \sigma_{x^l} y^l P =  bF^2 + 2 (a_{x^l}y^l + \sigma \theta) F  + 3 \theta^2 ,
\ee
\[ \theta = \lambda \sigma_{x^l}y^l. \]
\end{lem}

{\it Proof}: By Lemma \ref{lemthetaclose}, (\ref{aF}) holds. Differentiating it with respect to
$x^k$ yields
\[ a_{x^k} F^2 + 2 a F F_{x^k} -  \sigma_{x^l x^k}y^l F - \sigma_{x^l}y^l F_{x^k} + 2\theta_{x^k} P + 2 \theta P_{x^k} - \theta_{x^l x^k}y^l =0. \]

Set
\[ \Theta:= a F^2 - \sigma_{x^l}y^l F + 2\theta P - \theta_{x^l}y^l.  \]
By (\ref{aF}), $\Theta=0$.
Differentiating (\ref{aF}) with respect to $y^i$ yields
\be  \label{Theta_yi}
0 = \Theta_{\cdot i} = (2 a F  - \sigma_{x^l} y^l) F_{\cdot i} - F \sigma_{x^i} + 2 (\theta P_{\cdot i}  + P \theta_i  -  \theta_{x^i}) .
\ee
Here we use the fact that  $\theta$ is closed. Which implies $ \theta_{x^i} = \theta_{i x^l} y^l$.
Differentiating (\ref{Theta_yi}) with respect to $x^j$ yields
\be \label{Theta_ixj_1}
\begin{split}
0 = \Theta_{\cdot i x^j} = &(2 a_{x^j} F + 2 a F_{x^j} - \sigma_{x^j x^l} y^l ) F_{\cdot i} +(2 a F  - \sigma_{x^l} y^l) F_{\cdot i x^j} - F_{x^j} \sigma_{x^i} - F \sigma_{x^i x^j} \\
& + 2 ( \theta_{x^j}  P_{\cdot i} +   \theta P_{\cdot i x^j} + P_{x^j} \theta_i + P \theta_{i x^j}  - \theta_{x^i x^j} ). \\
\end{split}
\ee
By the assumption, $F$ is projectively flat, then by (\ref{PF-F}), (\ref{PF-P}) and (\ref{PF-Pkxl}) we have
\[  F_{x^j} = F P_{\cdot j} + P F_{\cdot j},  \]
\be
\begin{split}
 P_{\cdot i x^j} =&  \frac{1}{2}[P^2]_{\cdot i\cdot j} - F F_{\cdot j} {\bf K}_{\cdot i} - \frac{1}{2}{\bf K} [F^2]_{\cdot i\cdot j} - \frac{2}{3} F F_{\cdot i} {\bf K}_{\cdot j} -\frac{1}{3}F^2 {\bf K}_{\cdot i\cdot j} \\
 = & \frac{1}{2}[P^2]_{\cdot i\cdot j} - 2\theta F_{\cdot i\cdot j} - \frac{1}{2} \sigma [F^2]_{\cdot i\cdot j} - 2 \theta_{i} F_{\cdot j} - \theta_j F_{\cdot i}.
\end{split}
\ee
Plugging above two equations into (\ref{Theta_ixj_1}), we get
\be \label{Theta_ixj_2}
\begin{split}
0 =\ & \Theta_{\cdot i x^j} \\
=\  & (2 a_{x^j} F  + 2 a F P_{\cdot j} - \sigma_{x^j x^l}y^l - 2 \theta \theta_j ) F_{\cdot i} -  ( 8\theta \theta_i + 2 \sigma F \theta_i + P \sigma_{x^i}) F_{\cdot j}+ 2 \theta_{x^j}  P_{\cdot i} \\
&  + (2  P \theta_i -  F \sigma_{x^i}) P_{\cdot j} + 2 a P F_{\cdot i} F_{\cdot j} - 2 F \theta_i \theta_j - F \sigma_{x^i x^j} - 2 \theta_{x^i x^j}+ \theta [P^2]_{\cdot i \cdot j} \\
& - \sigma \theta [F^2]_{\cdot i \cdot j} - 4 \theta^2 F_{\cdot i \cdot j} + 2 P \theta_{i x^j}+(2 a F -\sigma_{x^l}y^l) F_{\cdot i x^j}. \\
\end{split}
\ee
Then we have
\be
\begin{split}
0 =\ & \Theta_{\cdot i x^j} - \Theta_{\cdot j x^i} \\
=\  &(2 a_{x^j} F  + 2 a F P_{\cdot j} - \sigma_{x^j x^l}y^l + 6 \theta \theta_j +  2 \sigma F \theta_j + P \sigma_{x^j}) F_{\cdot i} \\
& - (2 a_{x^i} F  + 2 a F P_{\cdot i} - \sigma_{x^i x^l}y^l + 6 \theta \theta_i +  2 \sigma F \theta_i + P \sigma_{x^i}) F_{\cdot j} \\
& + 2 \theta_{x^j}  P_{\cdot i} -  2 \theta_{x^i}  P_{\cdot j}+ (2  P \theta_i -  F \sigma_{x^i}) P_{\cdot j}- (2  P \theta_j -  F \sigma_{x^j}) P_{\cdot i}\\
\end{split}
\ee
Contracting above equation with $y^j$ yields
\be \label{Theta_zero}
\begin{split}
0 =\  &(2 a_{x^l}y^l  F  + 2 a F P - \sigma_{x^k x^l}y^k y^l + 6 \theta^2 +  2 \sigma \theta F + P \sigma_{x^l} y^l) F_{\cdot i} \\
& - (2 a_{x^i} F  + 2 a F P_{\cdot i} - \sigma_{x^i x^l}y^l + 6 \theta \theta_i +  2 \sigma F \theta_i + P \sigma_{x^i}) F  \\
& + 2 \theta_{x^l}y^l  P_{\cdot i} -  2 \theta_{x^i}  P + (2  P \theta_i  -  F \sigma_{x^i}) P - (2 \theta P -  F \sigma_{x^l}y^l) P_{\cdot i}. \\
\end{split}
\ee

Observing that
\[ 2 a_{x^l}y^l  F F_{\cdot i} - 2 a_{x^i} F^2 = - 2 F^3 [ \frac{ a_{x^l}y^l}{F}]_{\cdot i}, \ \ \  2 a F P F_{\cdot i} -2 a F^2 P_{\cdot i} = - 2 a F^3 [\frac{P}{F}]_{\cdot i}, \]

\[ - \sigma_{x^k x^l}y^k y^l F_{\cdot i} + F \sigma_{x^i x^l}y^l = \frac{1}{2} F^3 [\frac{\sigma_{x^k x^l}y^k y^l}{F^2}]_{\cdot i}, \ \ \ 6 \theta^2 F_{\cdot i} - 6 \theta F \theta_i = - 3 F^3 [\frac{\theta^2}{F^2}]_{\cdot i}, \]

\[ 2 \sigma \theta F F_{\cdot i} -  2 \sigma F^2 \theta_i = -2\sigma F^3 [\frac{\theta}{F}]_{\cdot i}, \ \ \  P \sigma_{x^l} y^l F_{\cdot i} + F \sigma_{x^l} y^l P_{\cdot i}- 2 F P \sigma_{x^i} = - 2 \sqrt{(FP)^3} [\frac{\sigma_{x^l} y^l}{\sqrt{PF}}]_{\cdot i}, \]

\[ 2 \theta_{x^l}y^l  P_{\cdot i} -  2 \theta_{x^i}  P = - P^3 [\frac{\theta_{x^l}y^l}{P^2}]_{\cdot i}, \ \ \  2  P^2 \theta_i - 2 \theta P P_{\cdot i} = 2 P^3 [\frac{\theta}{P}]_{\cdot i},\]

(\ref{Theta_zero}) can be written into
\[
\begin{split}
0 =\ & F^3 \big[ -2 \frac{ a_{x^l}y^l}{F} - 2 a \frac{P}{F} + \frac{1}{2}\frac{\sigma_{x^k x^l}y^k y^l}{F^2} -3 \frac{\theta^2}{F^2}-2\sigma\frac{\theta}{F} \big]_{\cdot i} \\
 & - 2 \sqrt{(FP)^3} [\frac{\sigma_{x^l} y^l}{\sqrt{PF}}]_{\cdot i}  + P^3 \big[ - \frac{\theta_{x^l}y^l}{P^2}  + 2 \frac{\theta}{P} \big]_{\cdot i}.
\end{split}
\]
It is equivalent to
\be \label{Theta_zero1}
\begin{split}
 &\big[ 2 \frac{ a_{x^l}y^l}{F} + 2 a \frac{P}{F} - \frac{1}{2}\frac{\sigma_{x^k x^l}y^k y^l}{F^2} + 3 \frac{\theta^2}{F^2} + 2\sigma\frac{\theta}{F} \big]_{\cdot i} \\
 =\ & - 2 (\frac{P}{F})^{\frac{3}{2}} [\frac{\sigma_{x^l} y^l}{\sqrt{PF}}]_{\cdot i}  + (\frac{P}{F})^3 \big[ - \frac{\theta_{x^l}y^l}{P^2}  + 2 \frac{\theta}{P} \big]_{\cdot i} \\
 =\ & - 2 (\frac{P}{F})^{\frac{3}{2}} [\frac{\sigma_{x^l} y^l}{P^2}(\frac{P}{F})^{\frac{3}{2}}]_{\cdot i}  + (\frac{P}{F})^3 \big[ \frac{\sigma_{x^l} y^l - a F^2}{P^2} \big]_{\cdot i}\\
 =\ & - [ \frac{\sigma_{x^l} y^l P}{F^2} ]_{\cdot i} + 2 a [\frac{P}{F}]_{\cdot i}.
\end{split}
\ee
Here (\ref{aF}) is used. Then there exists a scalar function $b = b(x)$ such that
\[ -b = 2 \frac{ a_{x^l}y^l}{F}  - \frac{1}{2}\frac{\sigma_{x^k x^l}y^k y^l}{F^2} + 3 \frac{\theta^2}{F^2} + 2\sigma\frac{\theta}{F} +\frac{\sigma_{x^l} y^l P}{F^2}.  \]
Thus
\be  \label{bF}
 -bF^2 - 2 (a_{x^l}y^l + \sigma \theta) F + \frac{1}{2}\sigma_{x^k x^l}y^k y^l  - 3 \theta^2   - \sigma_{x^l} y^l P =0.
\ee
By Lemma \ref{lemthetaclose}, (\ref{aF}) holds. Then $\sigma_{x^l} y^l \times (\ref{aF}) + 2 \theta \times (\ref{bF})$ yields
\be \label{proQ}
 (a \sigma_{x^l} y^l- 2 b \theta) F^2 - [ (\sigma_{x^l}y^l)^2 + 4 (a_{x^l}y^l + \sigma \theta)\theta  ]F    -\sigma_{x^k} y^k \theta_{x^l}y^l+ \theta \sigma_{x^k x^l}y^ky^l - 6 \theta^3=0.
\ee
If $(\sigma_{x^l}y^l)^2 + 4 (a_{x^l}y^l + \sigma \theta)\theta \neq 0$, then by Lemma \ref{lemtrivial} $F$ is a Randers metric.
If
\be  \label{coeF1}
(\sigma_{x^l}y^l)^2 + 4 (a_{x^l}y^l + \sigma \theta)\theta = 0,
\ee
 we divide the left proof
into two cases.

{\bf Case (i) $\sigma_{x^l}y^l =0$.} In this case, (\ref{bF}) becomes into
 \be
 - bF^2 - 2 (a_{x^l}y^l + \sigma \theta) F   - 3 \theta^2   =0.
\ee
Then by Lemma \ref{lemtrivial}, $F$ is a Randers metric.

{\bf Case (ii) $\sigma_{x^l}y^l \neq 0$.} In this case, by (\ref{coeF1}) there exists a scalar function
$\lambda  = \lambda(x)$ such that
\[ \theta = \lambda \sigma_{x^k} y^k. \]

\qed

{\bf Proof of Theorem \ref{thm2}}:
By Lemma \ref{lemprobF}, we only need to prove the theorem when $\sigma_{x^l}y^l = \sigma_{|0}\neq 0$. By the assumption $G^i = P y^i$, then (\ref{aF}) and (\ref{lembF}) can be rewritten into
\be \label{paF}
\theta_{|0} = aF^2 - \sigma_{|0} F
\ee
\be\label{psigma_v0v0}
 \frac{1}{2}\sigma_{| 0 |0} =bF^2 + 2 (a_{x^l}y^l + \sigma \theta) F + 3 \theta^2.
\ee
By a direct computation, we get
\be \label{thmpcase2}
\begin{split}
 \theta_{|0} =& \lambda_{|0} \sigma_{|m}y^m + \lambda \sigma_{|0|0} \\
                =&  2 \lambda b F^2 + 4 \lambda (a_{|0} +\sigma \theta)F + 6 \lambda \theta^2 + \lambda_{|0} \sigma_{|0}.
 \end{split}
\ee
Here the second equality is from (\ref{psigma_v0v0}). Substituting (\ref{thmpcase2}) back into (\ref{paF}) yields
\be \label{thmeqab}
 (a -2\lambda b) F^2 - [4\lambda(a_{|0} +\sigma \theta) + \sigma_{|0}]F - ( 6\lambda^3 \sigma_{|0} + \lambda_{|0}) \sigma_{|0} =0.
 \ee
If $4\lambda(a_{|0} +\sigma \theta) + \sigma_{|0} \neq 0$ or $6\lambda^3 \sigma_{|0} + \lambda_{|0}  \neq 0$, then by Lemma \ref{lemtrivial}, $F$ is a Randers metric. Otherwise
\be \label{thmpeq_last1} 4\lambda(a_{|0} +\sigma \theta) + \sigma_{|0} =0 \ee
and
 \be \label{thmpeq_last2} 6\lambda^3 \sigma_{|0} + \lambda_{|0}  =0. \ee
We claim that it is impossible if $F$ is not a Randers metric.
By replacing the function $f=f(x)$ with $a=a(x)$, (\ref{thmpeq_last1}), (\ref{thmpeq_last2}) and (\ref{paF}) are just
(\ref{lemlast1}), (\ref{lemlast2}) and (\ref{mainlemEq}) in Lemma \ref{lemlast}.
Then  Lemma \ref{lemlast}  is still true in this case (even in dimension two). Then
$\theta = 0$. This case is excluded by the assumption $\theta \neq 0$.
\qed

\section*{Acknowledgments}
The author would like to thank the referees  for their helpful comments. Research is supported by the National Science Foundation of China (11371209), National Science Foundation of Zhejiang Province (R18A010002, LY13A010013) and K.C. Wong Magna Fund in Ningbo University.

\end{document}